\documentclass[reqno]{amsart}
\usepackage[margin=1.25in]{geometry}
\usepackage[backref=page]{hyperref}
\usepackage{cite}
\usepackage{enumitem}

\theoremstyle{plain}
\newtheorem{theorem}{Theorem}[section]
\newtheorem{proposition}[theorem]{Proposition}
\newtheorem{corollary}[theorem]{Corollary}
\newtheorem{lemma}[theorem]{Lemma}

\theoremstyle{remark}
\newtheorem{remark}[theorem]{Remark}

\allowdisplaybreaks

\begin{document}

\author{Phuong Le}
\address{Phuong Le$^{1,2}$ (ORCID: 0000-0003-4724-7118)\newline
	$^1$Faculty of Economic Mathematics, University of Economics and Law, Ho Chi Minh City, Vietnam; \newline
	$^2$Vietnam National University, Ho Chi Minh City, Vietnam}
\email{phuongl@uel.edu.vn}


\title[The polyharmonic Hénon--Lane--Emden conjecture]{The polyharmonic Hénon--Lane--Emden conjecture}

\begin{abstract}
We study positive classical solutions of the polyharmonic Hénon--Lane--Emden system $(-\Delta)^m u=|x|^\alpha v^p$, $(-\Delta)^m v=|x|^\beta u^q$ in $\mathbb{R}^n$, where $m\ge2$, $p,q\ge1$, $pq>1$, and $\alpha,\beta\ge0$. We prove that there is no such solution whenever $\frac{n+\alpha}{p+1}+\frac{n+\beta}{q+1}>n-2m$, in every dimension and without any assumption on the behavior of the solutions at infinity. Combined with an existence result of Villavert, this proves the polyharmonic Hénon--Lane--Emden conjecture of Phan for $m\ge2$: for these exponents and weights, a positive classical solution exists if and only if $\frac{n+\alpha}{p+1}+\frac{n+\beta}{q+1}\le n-2m$. The proof has two parts. We first show that positivity forces $(-\Delta)^i u>0$ and $(-\Delta)^i v>0$ for $1\le i\le m-1$. For $m$ even this needs $pq>1$, and explicit radial examples show that it fails when $pq<1$ and when $p=q=1$. The system is then equivalent to a weighted integral system with a Riesz kernel of order $2m$. For such integral systems we prove the corresponding Liouville theorem for every Riesz order $\mu\in(1,n)$, adapting a recent argument for the second-order case that combines a localized virial identity with an interval comparison along source layers.
\end{abstract}

\maketitle

\section{Introduction}\label{sec_intro}

We are concerned with positive solutions of the system
\begin{equation}\label{eq_main}
(-\Delta)^m u=|x|^\alpha v^p,\qquad (-\Delta)^m v=|x|^\beta u^q\qquad\text{in } \mathbb{R}^n,
\end{equation}
where $m\ge2$ is an integer and $p,q>0$, $\alpha,\beta\in\mathbb{R}$. The scaling $u_R(x)=R^{\sigma_1}u(Rx)$, $v_R(x)=R^{\sigma_2}v(Rx)$ with
\begin{equation}\label{eq_sigma}
\sigma_1=\frac{2m+\alpha+p(2m+\beta)}{pq-1},\qquad \sigma_2=\frac{2m+\beta+q(2m+\alpha)}{pq-1}
\end{equation}
leaves \eqref{eq_main} invariant when $pq\neq1$, and it singles out the critical hyperbola
\begin{equation}\label{eq_hyperbola}
\frac{n+\alpha}{p+1}+\frac{n+\beta}{q+1}=n-2m.
\end{equation}
Phan \cite{Ph15} conjectured that, for $p,q\ge1$, $pq\neq1$, and $\alpha,\beta\ge0$, the system \eqref{eq_main} has no positive classical solution if and only if $(p,q)$ lies below the critical hyperbola. The nonexistence half, which is \cite[Conjecture A]{Ph15} and is also stated by Fazly \cite[Conjecture 1]{Fa13}, asserts that there is no positive classical solution when
\begin{equation}\label{eq_sub}
\frac{n+\alpha}{p+1}+\frac{n+\beta}{q+1}>n-2m.
\end{equation}
Fazly proved it for radial solutions in every dimension and for bounded solutions in dimension $n=2m+1$ \cite{Fa13}. Phan proved it for bounded solutions when, in addition, $\max\{2m(p+1),2m(q+1)\}\ge(n-2m-1)(pq-1)$, and without the boundedness assumption in a smaller range of exponents \cite{Ph15}. Liouville theorems for integral and polyharmonic systems have also been obtained by the method of scaling spheres \cite{DQ23,Le21}. Arthur and Yan obtained further results under a slow-decay assumption \cite{AY16}. More recently, Xu and Luo proved a super polyharmonic property and Liouville theorems for radial solutions subject to the growth condition $u,v=O(|x|^2)$ \cite{XL26}. These restrictions are absent from the result below. In higher dimensions the general unrestricted case has remained open.

The second-order case $m=1$ has a longer history. For $\alpha=\beta=0$ it is the Lane--Emden conjecture. Mitidieri introduced a Rellich identity for systems and obtained the first nonexistence regions \cite{Mi93,Mi96}. Serrin and Zou settled dimension three under a growth condition \cite{SZ96}. Poláčik, Quittner, and Souplet removed that condition \cite{PQS07}, and Souplet proved the conjecture in dimension four \cite{So09}. Phan formulated and partly proved the weighted version \cite{Ph12}, and Li and Zhang proved it in dimension three \cite{LZ19}. In each of these works the decisive step is an estimate for the boundary terms of a Rellich--Pohozaev identity, and the available Sobolev-type estimates control them only in low dimensions. A recent preprint \cite{OAI26} removes this restriction for $m=1$ and every $n$. Its main estimate does not bound the boundary term of the localized Pohozaev identity. Instead, it shows that this term is dominated by an interaction term of favorable sign, up to an error $\varepsilon E+C_\varepsilon$ relative to the localized energy $E$. The domination comes from an elementary comparison between pairs of intervals on parallel lines, applied to the level sets of the sources.

That argument uses the equation only through the Newton potential representation $u=\Gamma*(|x|^\alpha v^p)$, $v=\Gamma*(|x|^\beta u^q)$, together with the fact that the kernel is a power of $|x|$. Our first observation is that it goes through for a Riesz kernel of any order $\mu\in(1,n)$. Mei and Wei carried out this extension for the unweighted fractional system, that is, for Riesz kernels of order $2s\in(0,2)$ and $\alpha=\beta=0$ \cite{MW26}. When $2s\le1$ the kernel gradient is not locally integrable, and they combine the two equations before removing the diagonal. The present problem requires weights and orders $\mu=2m\ge4$. The natural route for $m\ge2$ is therefore to prove that positive solutions of \eqref{eq_main} are Riesz potentials of order $2m$.

This is where the polyharmonic problem differs from the second-order one. When $m=1$, a positive solution is superharmonic and the representation follows from the Riesz decomposition of superharmonic functions. When $m\ge2$, positivity of $u$ says nothing directly about the sign of $-\Delta u$, and the representation can fail. In fact, for $m$ even and $pq<1$, or for $p=q=1$, the system \eqref{eq_main} with $\alpha=\beta=0$ has positive radial solutions with $\Delta u>0$ everywhere, although \eqref{eq_sub} holds (Proposition \ref{prop_example} below). The property that has to be proved is the super polyharmonic property
\[
(-\Delta)^i u>0,\qquad (-\Delta)^i v>0\qquad\text{in } \mathbb{R}^n,\quad 1\le i\le m-1.
\]
For the scalar equation it goes back to Wei and Xu \cite{WX99}. For systems $(-\Delta)^ku_i=f_i(u_1,\dots,u_l)$ without weights, Chen and Li proved it and deduced the equivalence with the corresponding integral system when $2k<n$ \cite{CL13}, and Xu and Luo obtained it for radial solutions of weighted systems under a growth condition \cite{XL26}. The weights $|x|^\alpha$ and $|x|^\beta$ are not centered at the points where one averages, and we give a self-contained proof, without symmetry or growth assumptions, that also covers the endpoint cases $p=1$ or $q=1$. Spherical averages and a parity argument give the property at once except in one case: $m$ even and $(-\Delta)^{m-1}u$ negative somewhere. There the averages of $u$ and $v$ grow, but the parity count produces no contradiction. We show instead that the growth exponents feed each other and increase geometrically when $pq>1$, which is impossible for a finite function. The example above shows that $pq>1$ cannot be dropped here.

Our main result settles the nonexistence half of the conjecture.

\begin{theorem}\label{thm_main}
Let $m\ge2$ and $n\ge1$ be integers, let $p,q\ge1$ with $pq>1$, and let $\alpha,\beta\ge0$. If \eqref{eq_sub} holds, then \eqref{eq_main} has no solution $u,v\in C^{2m}(\mathbb{R}^n)$ with $u>0$ and $v>0$ in $\mathbb{R}^n$.
\end{theorem}

Together with the existence theorem of Villavert \cite{Vi14}, this settles the conjecture of Phan for $m\ge2$.

\begin{corollary}\label{cor_class}
Let $m\ge2$ and $n\ge1$ be integers, let $p,q\ge1$ with $pq\neq1$, and let $\alpha,\beta\ge0$. Then \eqref{eq_main} has a solution $u,v\in C^{2m}(\mathbb{R}^n)$ with $u>0$ and $v>0$ in $\mathbb{R}^n$ if and only if
\[
\frac{n+\alpha}{p+1}+\frac{n+\beta}{q+1}\le n-2m.
\]
In that case the solution can be chosen radial.
\end{corollary}

The radial solutions of \cite{Vi14} are constructed by shooting and are stated there to be of class $C^{2m}$ away from the origin. We check in Section \ref{sec_reduction} that they are classical solutions on all of $\mathbb{R}^n$.

No boundedness, radial symmetry, decay, or energy assumption is made in Theorem \ref{thm_main}. Taking $u=v$, $p=q$, and $\alpha=\beta$ gives the scalar statement: for $m\ge2$, $p>1$, $\alpha\ge0$, and $n\ge1$, the equation $(-\Delta)^m u=|x|^\alpha u^p$ has no positive classical solution in $\mathbb{R}^n$ when $p<(n+2m+2\alpha)/(n-2m)$, with the convention that the bound is $+\infty$ when $n\le2m$. For $n>2m$ this scalar statement was proved earlier, for $\alpha>-2m$ and $0<p$ below the same bound, by Dai and Qin with the method of scaling spheres \cite{DQ23}. When $n\le2m$, condition \eqref{eq_sub} holds automatically, and Theorem \ref{thm_main} says that there is no positive solution at all. The hypotheses $p,q\ge1$ and $\alpha,\beta\ge0$ are those of the conjecture. They enter the proof through the super polyharmonic property, which is the next statement, and $\alpha,\beta\ge0$ also makes the right sides of \eqref{eq_main} continuous at the origin.

\begin{theorem}\label{thm_super}
Let $m\ge2$ and $n\ge1$ be integers, let $p,q\ge1$ and $\alpha,\beta\ge0$, and let $u,v\in C^{2m}(\mathbb{R}^n)$ be positive solutions of \eqref{eq_main}. If $m$ is even, assume also $pq>1$. Then
\[
(-\Delta)^iu>0\quad\text{and}\quad(-\Delta)^iv>0\qquad\text{in }\mathbb{R}^n,\quad 1\le i\le m-1.
\]
\end{theorem}

\begin{proposition}\label{prop_example}
Let $m\ge2$ be even, $n\ge1$, and $\alpha=\beta=0$. Suppose that either $p,q>0$ with $pq<1$, or $p=q=1$. Then \eqref{eq_main} has a positive radial solution $u,v\in C^{2m}(\mathbb{R}^n)$ with $\Delta u\ge1$ and $\Delta v\ge1$ in $\mathbb{R}^n$. Moreover, $(p,q)$ satisfies \eqref{eq_sub}.
\end{proposition}

We do not know whether Theorem \ref{thm_super} holds when $pq>1$ and $\min\{p,q\}<1$. The proof below uses Jensen's inequality for the averages of $u^q$ and $v^p$, and this is the essential reason for the restriction $p,q\ge1$.

Once Theorem \ref{thm_super} is available, $u$ and $v$ solve an integral system, and Theorem \ref{thm_main} follows from the next result with $\mu=2m$. We state it for a general Riesz order because the proof is the same.

\begin{theorem}\label{thm_integral}
Let $n\ge2$, $1<\mu<n$, $c_0>0$, $p,q>0$, and $\alpha,\beta\in\mathbb{R}$. Suppose that
\begin{equation}\label{eq_sub_alpha}
\frac{n+\alpha}{p+1}+\frac{n+\beta}{q+1}>n-\mu.
\end{equation}
Then there are no functions $u,v\in C(\mathbb{R}^n)$, positive in $\mathbb{R}^n$, such that for every $x\in\mathbb{R}^n$ the integrals below are finite and
\begin{equation}\label{eq_int}
u(x)=c_0\int_{\mathbb{R}^n}\frac{|y|^\alpha v(y)^p}{|x-y|^{n-\mu}}\,dy,\qquad v(x)=c_0\int_{\mathbb{R}^n}\frac{|y|^\beta u(y)^q}{|x-y|^{n-\mu}}\,dy.
\end{equation}
\end{theorem}

For $\mu=2$ and $n\ge3$, Theorem \ref{thm_integral} is the integral form of the main theorem of \cite{OAI26}, and its proof is the proof given there. We give in full the steps in which the order of the kernel enters. The restriction $\mu>1$ is used only to differentiate the potentials in the virial identity, where $\nabla|z|^{\mu-n}$ must be locally integrable. For $\alpha=\beta=0$ the system \eqref{eq_int} is the integral system that Chen and Li related to the Lane--Emden conjecture \cite{CL09}, and for $\alpha=\beta=0$ and $\mu\in(1,2)$ Theorem \ref{thm_integral} is contained in \cite{MW26}, after normalizing $c_0$ by scaling. We do not know whether the method of \cite{MW26} for $\mu\le1$ extends to the weighted case.

Our integral-system argument builds on the recent preprint \cite{OAI26}. For the Riesz kernel $|z|^{\mu-n}$, we establish the interval comparison and the localized estimates in Section \ref{sec_integral}, followed by the weighted pressure estimate in Appendix \ref{app_pressure}. We also give the virial identity with the required kernel-gradient integrability arguments. By contrast, the super polyharmonic property, the reduction to \eqref{eq_int}, and Proposition \ref{prop_example} do not depend on \cite{OAI26}.

Theorem \ref{thm_super} is proved in Section \ref{sec_super}, the reduction to the integral system and the proofs of Theorem \ref{thm_main} and Corollary \ref{cor_class} are in Section \ref{sec_reduction}, Theorem \ref{thm_integral} is proved in Section \ref{sec_integral}, and Proposition \ref{prop_example} in Section \ref{sec_example}.

\section{The super polyharmonic property}\label{sec_super}

Throughout this section $u,v$ are as in Theorem \ref{thm_super}, and we write
\[
u_i=(-\Delta)^iu,\qquad v_i=(-\Delta)^iv\qquad(0\le i\le m-1),\qquad u_m=|x|^\alpha v^p,\qquad v_m=|x|^\beta u^q.
\]
Thus $u_i,v_i\in C^{2m-2i}(\mathbb{R}^n)$ and $-\Delta u_i=u_{i+1}$, $-\Delta v_i=v_{i+1}$ for $0\le i\le m-1$. The functions $u_m$ and $v_m$ are continuous and positive away from the origin.

Fix $x_0\in\mathbb{R}^n$. For a continuous function $w$ we write
\[
\bar w(r)=\frac{1}{|\partial B_r(x_0)|}\int_{\partial B_r(x_0)}w\,dS\qquad(r>0),\qquad \bar w(0)=w(x_0),
\]
with the counting measure on $\partial B_r(x_0)=\{x_0\pm r\}$ when $n=1$. If $w\in C^2(\mathbb{R}^n)$, the divergence theorem gives $r^{n-1}\bar w'(r)=|\partial B_1|^{-1}\int_{B_r(x_0)}\Delta w\,dx$. Hence, for $0\le i\le m-1$,
\begin{equation}\label{eq_avg}
r^{n-1}\bar u_i'(r)=-\int_0^rt^{n-1}\bar u_{i+1}(t)\,dt,\qquad \bar u_i'(0)=0,
\end{equation}
and the same holds for $v$. Two further facts are used repeatedly. First, if $r\ge R\ge\max\{2,2|x_0|\}$ and $x\in\partial B_r(x_0)$, then $|x|\ge r/2\ge1$. Since $\alpha,\beta\ge0$ and $p,q\ge1$, Jensen's inequality gives
\begin{equation}\label{eq_jensen}
\bar u_m(r)\ge\bar v(r)^p,\qquad \bar v_m(r)\ge\bar u(r)^q\qquad(r\ge\max\{2,2|x_0|\}).
\end{equation}
Second, $\bar u_m\ge0$ and $\bar v_m\ge0$ everywhere.

We first describe how a negative value propagates down a chain of averages.

\begin{lemma}\label{lem_chain}
Let $\ell\ge2$, let $W_0,\dots,W_{\ell-1}\in C^2([0,\infty))$ and $W_\ell\in C([0,\infty))$ satisfy
\[
r^{n-1}W_i'(r)=-\int_0^rt^{n-1}W_{i+1}(t)\,dt\qquad(0\le i\le\ell-1),
\]
and suppose $W_\ell\ge0$ and $W_{\ell-1}(r_1)<0$ for some $r_1\ge0$. Put $\epsilon_i=(-1)^{\ell-i}$. Then there is $R_*$ such that, for $0\le i\le\ell-1$ and $r\ge R_*$,
\[
\epsilon_iW_i(r)>0,\qquad \epsilon_iW_i'(r)\ge0,
\]
and $\epsilon_iW_i(r)\to\infty$ as $r\to\infty$ for $0\le i\le\ell-2$.
\end{lemma}

\begin{proof}
Since $W_\ell\ge0$, the function $W_{\ell-1}$ is nonincreasing, so $\epsilon_{\ell-1}W_{\ell-1}=-W_{\ell-1}\ge c:=-W_{\ell-1}(r_1)>0$ on $[r_1,\infty)$ and $\epsilon_{\ell-1}W_{\ell-1}'\ge0$ everywhere. Suppose that for some $i\le\ell-2$ we have $\epsilon_{i+1}W_{i+1}\ge c'>0$ on $[R,\infty)$. Since $\epsilon_i=-\epsilon_{i+1}$,
\[
\epsilon_ir^{n-1}W_i'(r)=\int_0^rt^{n-1}\epsilon_{i+1}W_{i+1}(t)\,dt\ge-\int_0^Rt^{n-1}|W_{i+1}|\,dt+\frac{c'}{n}(r^n-R^n),
\]
so $\epsilon_iW_i'(r)\ge c'r/(2n)$ for large $r$, and $\epsilon_iW_i(r)\to\infty$. Descending from $i=\ell-2$ to $i=0$ proves the lemma.
\end{proof}

The second lemma converts a lower bound for the top of a chain into a lower bound for its bottom on a fixed window $[R,2R]$.

\begin{lemma}\label{lem_window}
Let $\ell\ge1$, $R>0$, and let $W_0,\dots,W_\ell$ be as in the first display of Lemma \ref{lem_chain}. Let $\epsilon_i\in\{-1,1\}$ with $\epsilon_i=-\epsilon_{i+1}$, and suppose $\epsilon_iW_i(R)\ge0$ and $\epsilon_iW_i'(R)\ge0$ for $0\le i\le\ell-1$. If $\epsilon_\ell W_\ell(\tau)\ge C(\tau-R)^s$ on $[R,2R]$ for some $C>0$ and $s\ge0$, then
\[
\epsilon_0W_0(r)\ge C\,2^{\ell(1-n)}(s+2\ell)^{-2\ell}(r-R)^{s+2\ell}\qquad(R\le r\le2R).
\]
\end{lemma}

\begin{proof}
For $0\le i\le\ell-1$ and $r\ge R$, integrating the relation for $W_i$ from $R$ gives
\[
\epsilon_iW_i(r)=\epsilon_iW_i(R)+\epsilon_iR^{n-1}W_i'(R)\int_R^rt^{1-n}\,dt+\int_R^rt^{1-n}\int_R^t\tau^{n-1}\epsilon_{i+1}W_{i+1}(\tau)\,d\tau\,dt.
\]
The first two terms are nonnegative. On $[R,2R]$ we have $\tau^{n-1}t^{1-n}\ge2^{1-n}$ for $\tau\le t$, so a lower bound $\epsilon_{i+1}W_{i+1}(\tau)\ge C'(\tau-R)^{s'}$ on $[R,2R]$ gives
\[
\epsilon_iW_i(r)\ge\frac{2^{1-n}C'}{(s'+1)(s'+2)}(r-R)^{s'+2}\qquad(R\le r\le2R).
\]
Applying this $\ell$ times, from $i=\ell-1$ down to $i=0$, and using $(s+2j-1)(s+2j)\le(s+2\ell)^2$ for $j\le\ell$ proves the claim.
\end{proof}

\begin{proof}[Proof of Theorem \ref{thm_super}]
\textit{Step 1: nonnegativity.} Suppose that some $u_i$ or $v_i$ with $1\le i\le m-1$ takes a negative value, and let $k$ be the largest such index. Interchanging $(u,p,\alpha)$ and $(v,q,\beta)$ if necessary, we may assume $u_k(x_0)<0$ for some $x_0$. By maximality, $u_i\ge0$ and $v_i\ge0$ in $\mathbb{R}^n$ for $k<i\le m-1$. In particular $\bar u_{k+1}\ge0$, whether $k<m-1$ or $k=m-1$.

Apply Lemma \ref{lem_chain} to $W_i=\bar u_i$, $0\le i\le k$, with $\ell=k+1\ge2$ and $r_1=0$. It gives $(-1)^{k+1}\bar u(r)\to\infty$. If $k$ is even, then $\bar u(r)\to-\infty$, which is impossible because $u>0$.

Let $k$ be odd. Then $\bar u(r)\to\infty$, and \eqref{eq_jensen} gives $\bar v_m(r)\ge1$ for $r$ large. By \eqref{eq_avg}, $r^{n-1}\bar v_{m-1}'(r)\le-c(r^n-R^n)$ for $r\ge R$ if $R$ is large, so $\bar v_{m-1}(r)\to-\infty$. If $k<m-1$, this contradicts $v_{m-1}\ge0$.

It remains to treat $k=m-1$ with $m$ even. Here we use $pq>1$. Lemma \ref{lem_chain}, applied to $(\bar u_0,\dots,\bar u_m)$ with $\ell=m$ and $r_1=0$, and to $(\bar v_0,\dots,\bar v_m)$ with $\ell=m$ and $r_1$ chosen with $\bar v_{m-1}(r_1)<0$, gives $R_*$ such that, with $\epsilon_i=(-1)^{m-i}$,
\[
\epsilon_i\bar u_i,\ \epsilon_i\bar u_i',\ \epsilon_i\bar v_i,\ \epsilon_i\bar v_i'\ge0\quad\text{on }[R_*,\infty),\qquad 0\le i\le m-1.
\]
Put $c_1=-u_{m-1}(x_0)>0$. Since $\bar u_{m-1}$ is nonincreasing, we have $-\bar u_{m-1}\ge c_1$ on $[0,\infty)$.
Choose $R\ge\max\{R_*,2,2|x_0|\}$, to be fixed below. Applying Lemma \ref{lem_window} to $\bar u_0,\dots,\bar u_{m-1}$, with $\ell=m-1$ and $s=0$, gives
\[
\begin{aligned}
\bar u(r)&\ge C_0(r-R)^{s_0}\quad (R\le r\le 2R),\qquad s_0=2m-2,\\
C_0&=c_1\,2^{(m-1)(1-n)}(2m-2)^{-(2m-2)}.
\end{aligned}
\]
Suppose that $\bar u(r)\ge C_j(r-R)^{s_j}$ on $[R,2R]$. By \eqref{eq_jensen}, $\bar v_m\ge C_j^q(r-R)^{qs_j}$ there. Lemma \ref{lem_window} with $\ell=m$ gives a lower bound for $\bar v$, then \eqref{eq_jensen} gives one for $\bar u_m$, and Lemma \ref{lem_window} once more gives $\bar u(r)\ge C_{j+1}(r-R)^{s_{j+1}}$ on $[R,2R]$, where
\[
s_{j+1}=p(qs_j+2m)+2m,\qquad C_{j+1}=C_j^{pq}\,2^{m(1-n)(p+1)}s_{j+1}^{-2m(p+1)}.
\]
Here we used $qs_j+2m\le s_{j+1}$, which holds because $p\ge1$. The constants $C_0$, $s_0$, and the recursion do not depend on $R$.

Set $L_j=\log(C_jR^{s_j})$, so that $\bar u(2R)\ge e^{L_j}$. Since $s_{j+1}-pqs_j=2m(p+1)$ and $R\ge1$,
\[
L_{j+1}\ge pqL_j-b-2m(p+1)\log s_{j+1},\qquad b=m(n-1)(p+1)\log2.
\]
From $s_{j+1}\le pq\,s_j+2m(p+1)$ we get $s_j\le S(pq)^j$ with $S=s_0+2m(p+1)/(pq-1)$. Hence
\[
L_{j+1}\ge pqL_j-a(j+1)-b',\qquad a=2m(p+1)\log(pq),\quad b'=b+2m(p+1)\log S.
\]
Let $\theta=(1+pq)/2\in(1,pq)$ and $M=\max\{a+b',a/\log\theta\}/(pq-\theta)$. If $L_j\ge M\theta^j$, then, since $\theta^j\ge1+j\log\theta$,
\[
L_{j+1}\ge\theta^{j+1}M+(pq-\theta)M\theta^j-a(j+1)-b'\ge\theta^{j+1}M.
\]
Now $L_0=\log C_0+(2m-2)\log R$ and $2m-2>0$. Fix $R$ so large that $L_0\ge M$. Then $L_j\ge M\theta^j$ for every $j$, and $\bar u(2R)=\infty$, which is absurd. This proves $u_i\ge0$ and $v_i\ge0$ for $1\le i\le m-1$.

\medskip
\noindent\textit{Step 2: strict positivity.} Let $1\le i\le m-1$. The function $u_i$ is nonnegative and $-\Delta u_i=u_{i+1}\ge0$. If $u_i(x_1)=0$ at some point, the strong minimum principle for superharmonic functions \cite{AG01} gives $u_i\equiv0$, hence $u_{i+1}\equiv0$, and descending in this way we reach $|x|^\alpha v^p\equiv0$, which is false. The same argument applies to $v_i$.
\end{proof}

\begin{remark}\label{rem_parity}
The case $k=m-1$ with $m$ even is the only place where $pq>1$ is used. For $m$ odd, Theorem \ref{thm_super} holds for all $p,q\ge1$. For $m$ even it fails for $p=q=1$, by Proposition \ref{prop_example}.
\end{remark}

\begin{remark}\label{rem_nonneg}
Fazly states the conjecture for nonnegative solutions \cite{Fa13}. Step 1 uses only $u,v\ge0$. In Step 2, the same minimum principle applied to $u$ and $v$ shows that a nonnegative solution is either positive or identically zero. Hence Theorems \ref{thm_super} and \ref{thm_main} hold for nonnegative solutions, with the conclusion $u\equiv v\equiv0$ in Theorem \ref{thm_main}.
\end{remark}

\section{Reduction to the integral system}\label{sec_reduction}

We need the Riesz decomposition of a nonnegative superharmonic function on $\mathbb{R}^n$ in the following elementary form. Let $\Gamma(z)=|z|^{2-n}/((n-2)|\partial B_1|)$ for $n\ge3$.

\begin{lemma}\label{lem_riesz}
Let $w\in C^2(\mathbb{R}^n)$ satisfy $w\ge0$ and $-\Delta w=h\ge0$ with $h$ continuous.
\begin{enumerate}[label=\textup{(\roman*)}]
\item If $n\le2$, then $h\equiv0$.
\item If $n\ge3$, then $\Gamma*h$ is finite everywhere and $w=c+\Gamma*h$ for a constant $c\ge0$.
\end{enumerate}
\end{lemma}

\begin{proof}
(i) Average about a point $x_0$ with $h(x_0)>0$. Then $(r^{n-1}\bar w')'=-r^{n-1}\bar h\le0$, and $r^{n-1}\bar w'(r)\le-c<0$ for $r\ge r_1$. For $n=2$ this gives $\bar w(r)\le\bar w(r_1)-c\log(r/r_1)$, and for $n=1$ it gives $\bar w(r)\le\bar w(r_1)-c(r-r_1)$. Both contradict $w\ge0$.

(ii) For every $x$ and $R>0$, Green's representation in $B_R(x)$ gives
\[
w(x)-\bar w(R)=\int_{B_R(x)}\big(\Gamma(x-y)-\Gamma(R)\big)h(y)\,dy,
\]
where $\bar w$ is the average about $x$. Since $\bar w(R)\ge0$, monotone convergence as $R\to\infty$ gives $\Gamma*h\le w$. The rest of the argument is that of \cite[Lemma 2.2]{OAI26}: $\Gamma*h$ is locally integrable and solves $-\Delta(\Gamma*h)=h$ in the sense of distributions, so $w-\Gamma*h$ is a nonnegative harmonic function on $\mathbb{R}^n$, hence a nonnegative constant $c$ by Liouville's theorem for positive harmonic functions \cite{AG01}. A priori this gives $w=c+\Gamma*h$ almost everywhere. Since $\Gamma*h$ is superharmonic, it equals at every point the limit of its averages over small balls, and so does the continuous function $w-c$. Hence $w=c+\Gamma*h$ everywhere.
\end{proof}

For integers $j\ge1$ with $2j<n$ let $\Gamma_{2j}=\Gamma*\cdots*\Gamma$ be the $j$-fold convolution. It is the Riesz kernel $\Gamma_{2j}(z)=\gamma_{n,2j}|z|^{2j-n}$ with a constant $\gamma_{n,2j}>0$ \cite[Chapter V, \S1]{St70}, and $\Gamma*(\Gamma_{2j-2}*h)=\Gamma_{2j}*h$ for every nonnegative $h$, by Tonelli's theorem.

\begin{proposition}\label{prop_reduction}
Let $m$, $n$, $p$, $q$, $\alpha$, $\beta$ be as in Theorem \ref{thm_main}, and let $u,v\in C^{2m}(\mathbb{R}^n)$ be positive solutions of \eqref{eq_main}. Then $n>2m$, and for every $x\in\mathbb{R}^n$
\begin{equation}\label{eq_repr}
u(x)=\int_{\mathbb{R}^n}\Gamma_{2m}(x-y)|y|^\alpha v(y)^p\,dy,\qquad v(x)=\int_{\mathbb{R}^n}\Gamma_{2m}(x-y)|y|^\beta u(y)^q\,dy.
\end{equation}
\end{proposition}

\begin{proof}
By Theorem \ref{thm_super}, $u_i>0$ and $v_i>0$ for $0\le i\le m$ away from the origin, with the notation of Section \ref{sec_super}. Put $g=u_m=|x|^\alpha v^p$.

If $n\le2$, Lemma \ref{lem_riesz}(i) applied to $w=u_{m-1}$ gives $g\equiv0$, which is false. Let $n\ge3$. By Lemma \ref{lem_riesz}(ii), $u_{m-1}=c_{m-1}+\Gamma*g$ with $c_{m-1}\ge0$. We claim that for $1\le j\le m$ with $2j<n$,
\begin{equation}\label{eq_induct}
u_{m-j}=c_{m-j}+\Gamma_{2j}*g,\qquad c_{m-j}\ge0.
\end{equation}
Suppose \eqref{eq_induct} holds for $j-1\ge1$. Lemma \ref{lem_riesz}(ii) applied to $w=u_{m-j}$ and $h=u_{m-j+1}$ gives $\Gamma*u_{m-j+1}\le u_{m-j}<\infty$. Since $\Gamma*1=\infty$, this forces $c_{m-j+1}=0$, and then $\Gamma*u_{m-j+1}=\Gamma_{2j}*g$, which is \eqref{eq_induct} for $j$.

Suppose $2j\ge n$ for some $j\le m$, and take the smallest such $j$. Then $j\ge2$ because $n\ge3$. As above, $\Gamma*(\Gamma_{2j-2}*g)\le u_{m-j}<\infty$. On the other hand, for $|y|\ge1$ we have $|y-z|\le2|y|$ when $|z|\le1$, so
\[
\Gamma_{2j-2}*g(y)\ge\gamma_{n,2j-2}(2|y|)^{2j-2-n}\int_{B_1}g=:c|y|^{2j-2-n},\qquad c>0.
\]
The integral of $|x-y|^{2-n}|y|^{2j-2-n}$ over $\{|y|\ge2|x|+1\}$ diverges, because the integrand is comparable to $|y|^{2j-2n}$ and $2j-2n\ge-n$. This contradiction shows $n>2m$. Then \eqref{eq_induct} holds for $j=m$, and $u=c_u+\Gamma_{2m}*g$ with $c_u\ge0$. These identities hold at every point, since each $\Gamma_{2j}*g=\Gamma*(\Gamma_{2j-2}*g)$ is superharmonic and the argument at the end of the proof of Lemma \ref{lem_riesz} applies. In the same way $v=c_v+\Gamma_{2m}*f$ with $f=|x|^\beta u^q$.

It remains to show $c_u=c_v=0$. Since $\Gamma_{2m}(x-y)\ge\gamma_{n,2m}(2R)^{2m-n}$ for $x,y\in B_R$, we have
\[
u\ge\gamma_{n,2m}(2R)^{2m-n}\int_{B_R}g,\qquad v\ge\gamma_{n,2m}(2R)^{2m-n}\int_{B_R}f\qquad\text{in }B_R.
\]
Exactly as in \cite[Proposition 2.4]{OAI26}, the case $R=1$ gives $\int_{B_1}f\ge c\big(\int_{B_1}g\big)^q$ and $\int_{B_1}g\ge c\big(\int_{B_1}f\big)^p$ with $c$ depending only on $n,m,p,q,\alpha,\beta$. As $pq>1$, this bounds $\int_{B_1}f$ by a constant independent of the solution. The pair $u_R(x)=R^{\sigma_1}u(Rx)$, $v_R(x)=R^{\sigma_2}v(Rx)$ is again a positive $C^{2m}$ solution of \eqref{eq_main}, and changing variables yields $\int_{B_R}f\le CR^{n-2m-\sigma_2}$ for all $R>0$. If $c_u>0$, then $\int_{B_R}f\ge c_u^q\int_{B_R}|x|^\beta\,dx=cR^{n+\beta}$, and $n+\beta>n-2m-\sigma_2$, a contradiction for large $R$. Thus $c_u=0$, and likewise $c_v=0$.
\end{proof}

\begin{proof}[Proof of Theorem \ref{thm_main}]
Let $u,v$ be as in Theorem \ref{thm_main}. By Proposition \ref{prop_reduction}, $n>2m\ge4$ and $(u,v)$ satisfies \eqref{eq_int} with $\mu=2m\in(1,n)$ and $c_0=\gamma_{n,2m}$. Condition \eqref{eq_sub} is \eqref{eq_sub_alpha} with $\mu=2m$. Theorem \ref{thm_integral} shows that no such pair exists.
\end{proof}

\begin{proof}[Proof of Corollary \ref{cor_class}]
Since $p,q\ge1$ and $pq\neq1$, we have $pq>1$, and Theorem \ref{thm_main} shows that the condition is necessary. Conversely, assume the condition. Its left side is positive, so $n>2m$. We apply \cite[Corollary 2]{Vi14}, which rests on \cite[Theorems 1 and 3]{Vi14}: a shooting argument, and the nonexistence of radial solutions of the Navier problem on balls, proved there by a Pohozaev identity. With its order taken equal to $2m$, it concerns the system $(-\Delta)^mu=|x|^{-s_1}v^{q}$, $(-\Delta)^mv=|x|^{-s_2}u^{p}$ with $s_1,s_2<2$, where the weights are denoted $\sigma_1,\sigma_2$ in \cite{Vi14}. We take $s_1=-\alpha$, $s_2=-\beta$, and interchange the names of $p$ and $q$. The exponent condition there then becomes ours, and the hypothesis $2m<n$ of \cite{Vi14} holds. It gives a radial pair $u,v$, positive in $\mathbb{R}^n$ and of class $C^{2m}(\mathbb{R}^n\setminus\{0\})$, solving \eqref{eq_main} away from the origin.

It remains to check the origin. In \cite{Vi14} the pair is obtained from the second-order system satisfied by $w_i=(-\Delta)^iu$ and $z_i=(-\Delta)^iv$, $0\le i\le m-1$, with initial data $w_i(0),z_i(0)>0$ and $w_i'(0)=z_i'(0)=0$. Hence each component $W$ satisfies
\[
W(r)=W(0)-\int_0^rt^{1-n}\int_0^t\tau^{n-1}F(\tau)\,d\tau\,dt,
\]
where $F$ is the next component, or $|x|^\alpha v^p$, or $|x|^\beta u^q$. All of these are continuous on $[0,\infty)$ because $\alpha,\beta\ge0$. As noted in Section \ref{sec_example}, $W(|x|)$ is then of class $C^1$ and solves $-\Delta W=F$ in the sense of distributions on $\mathbb{R}^n$. In particular $u,v\in C^1(\mathbb{R}^n)$, so $|x|^\alpha v^p$ and $|x|^\beta u^q$ are locally Hölder continuous. Applying Schauder estimates to $w_{m-1},w_{m-2},\dots,w_0$ in turn, and likewise to the $z_i$, gives $u,v\in C^{2m}(\mathbb{R}^n)$.
\end{proof}

\begin{remark}\label{rem_linear}
For $m$ odd, the system \eqref{eq_main} with $p=q=1$ has no positive classical solution, for any $n\ge1$ and $\alpha,\beta\ge0$. Indeed, by Remark \ref{rem_parity}, Theorem \ref{thm_super} applies. The proof of Proposition \ref{prop_reduction}, up to the determination of $c_u$ and $c_v$, does not use $pq>1$. It gives $n>2m$, $u\ge\Gamma_{2m}*(|x|^\alpha v^p)$, and $v\ge\Gamma_{2m}*(|x|^\beta u^q)$. The argument of Step 1 in Section \ref{sec_integral} uses only these inequalities and forces $pq>1$, which is false. By Proposition \ref{prop_example}, for $m$ even and $\alpha=\beta=0$ such solutions exist. So, at least without weights, the parity of $m$ decides the linear case.
\end{remark}

\section{The integral system}\label{sec_integral}

In this section $n\ge2$, $1<\mu<n$, and $(u,v)$ is a positive continuous solution of \eqref{eq_int}. We write
\[
K(z)=c_0|z|^{\mu-n},\qquad f=|x|^\beta u^q,\qquad g=|x|^\alpha v^p,\qquad a=\frac1{q+1},\qquad b=\frac1{p+1},
\]
so that $u=K*g$ and $v=K*f$ at every point. We follow \cite{OAI26}. We follow the strategy of \cite{OAI26} but give the argument for the general Riesz kernel, including the virial identity and the pressure estimate in Appendix \ref{app_pressure}. The second-order kernel there has $n-2$ in place of $n-\mu$, and a statement of \cite{OAI26} quoted below is meant with this replacement. The weights $\alpha,\beta$ and the scaling exponents $\sigma_1,\sigma_2$ are denoted $A,B$ and $\alpha,\beta$ in \cite{OAI26}. To avoid a further clash, we write $m_0=n+3$ for the exponent denoted $m$ in \cite[(3.2)]{OAI26}. Write $k=n(m_0+1)$, $N=k+n+2$, and take the radial cutoffs
\[
d(x)=\min\{1,(2-|x|)_+\},\qquad
\eta(x)=\big[(1-(|x|-1)_+^3)_+\big]^N,\qquad
w(x)=-x\cdot\nabla\eta(x),\qquad X(x)=x\eta(x).
\]
These coincide with the cutoffs of \cite[(3.2)]{OAI26}; their differential bounds are collected in Appendix \ref{app_pressure}.

\medskip
\noindent\textit{Step 1: origin and exponents.} Since $v(0)>0$ and $v$ is continuous, $v\ge v(0)/2$ in some ball $B_\rho$. Finiteness of $u(0)$ then gives $\int_{B_\rho}|y|^{\mu-n+\alpha}\,dy<\infty$, that is, $\alpha>-\mu$. Likewise $\beta>-\mu$. In particular $f,g\in L^1_{\rm loc}(\mathbb{R}^n)$. For $x,y\in B_R$ we have $K(x-y)\ge c_0(2R)^{\mu-n}$. With $m_u(R)=\min_{\overline{B_R}}u$ and $m_v(R)=\min_{\overline{B_R}}v$, this gives
\[
m_u(R)\ge cR^{\mu+\alpha}m_v(R)^p,\qquad m_v(R)\ge cR^{\mu+\beta}m_u(R)^q.
\]
Hence $m_u(R)^{1-pq}\ge cR^{\mu+\alpha+p(\mu+\beta)}$, with a positive exponent. If $pq\le1$, the left side stays bounded as $R\to\infty$, so $pq>1$. We set
\[
\sigma_1=\frac{\mu+\alpha+p(\mu+\beta)}{pq-1},\qquad\sigma_2=\frac{\mu+\beta+q(\mu+\alpha)}{pq-1}.
\]
Then $p\sigma_2=\sigma_1+\mu+\alpha$ and $q\sigma_1=\sigma_2+\mu+\beta$, and the pair $u_R(x)=R^{\sigma_1}u(Rx)$, $v_R(x)=R^{\sigma_2}v(Rx)$ is again a positive continuous solution of \eqref{eq_int}. A direct computation gives
\begin{equation}\label{eq_gap}
\gamma_\mu:=a(n+\beta)+b(n+\alpha)-(n-\mu)=(1-a-b)(\sigma_1+\sigma_2+\mu-n),\qquad 1-a-b=\frac{pq-1}{(p+1)(q+1)},
\end{equation}
so \eqref{eq_sub_alpha} means $\gamma_\mu>0$ and implies $\sigma_1+\sigma_2+\mu-n>0$. Finally,
\begin{equation}\label{eq_energy_scaling}
\int_{B_1}(f_Ru_R+g_Rv_R)=R^{\sigma_1+\sigma_2+\mu-n}\int_{B_R}(fu+gv).
\end{equation}

\medskip
\noindent\textit{Step 2: universal bounds.} We make the mass estimate explicit, since it is used repeatedly below. Put $F=\int_{B_1}f$ and $G=\int_{B_1}g$. Both numbers are positive and finite. The lower bound $K(x-y)\ge c_0 2^{\mu-n}$ for $x,y\in B_1$ gives $u\ge cG$ and $v\ge cF$ on $B_1$, with $c>0$ depending only on the kernel. Consequently,
\[
F=\int_{B_1}|x|^\beta u^q\,dx\ge c_1G^q,\qquad
G=\int_{B_1}|x|^\alpha v^p\,dx\ge c_2F^p.
\]
Since $pq>1$, we have $F^{pq-1}\le (c_1c_2^q)^{-1}$ and, symmetrically, $G^{pq-1}\le(c_2c_1^p)^{-1}$. Thus $F+G\le C$ independently of the solution. Apply this bound to $(u_R,v_R)$ from Step 1. The identities $f_R(x)=R^{\sigma_2+\mu}f(Rx)$ and $g_R(x)=R^{\sigma_1+\mu}g(Rx)$ give, for every $R>0$,
\begin{equation}\label{eq_mass}
\int_{B_R}f\le CR^{n-\mu-\sigma_2},\qquad\int_{B_R}g\le CR^{n-\mu-\sigma_1},
\end{equation}
\begin{equation}\label{eq_tail}
\int_{B_4}(f+g+u+v)\le C,\qquad\int_{|x-y|\ge\ell}K(x-y)\big(f(y)+g(y)\big)\,dy\le C\ell^{\mu-n}\quad(x\in B_2,\ 0<\ell\le1).
\end{equation}
In the dyadic tail estimate, the annulus $2^j\le|y|<2^{j+1}$ contributes at most $C2^{j(\mu-n)}2^{(j+1)(n-\mu-\sigma_2)}\le C2^{-j\sigma_2}$, and the series converges because $\sigma_2>0$. The same dyadic estimate holds for $x\in B_4$, because $|x-y|\ge2^{j-1}$ there for $j\ge3$. Together with $K\in L^1_{\rm loc}$ and Tonelli's theorem, it gives the bound for $\int_{B_4}(u+v)$. For the second estimate in \eqref{eq_tail}, the part of the integral over $y\in B_8$ is at most $c_0\ell^{\mu-n}\int_{B_8}(f+g)\le C\ell^{\mu-n}$ by \eqref{eq_mass}, and the part over $y\notin B_8$ is at most $C\le C\ell^{\mu-n}$ by the dyadic estimate, since $\ell^{\mu-n}\ge1$ for $0<\ell\le1$.

\medskip
\noindent\textit{Step 3: differentiability.} Here $\mu>1$ is used. The function $|\nabla K(z)|=(n-\mu)c_0|z|^{\mu-n-1}$ is locally integrable, and $g$ is bounded near every $x_0\neq0$. Split $g$ into a bounded part with compact support near $x_0$ and a remainder, as in the last paragraph of the proof of \cite[Lemma 2.2]{OAI26}. Since $K\in W^{1,1}_{\rm loc}$, the convolution of $K$ with the first part is $C^1$ with gradient given by $\nabla K$, and the remainder is handled by dominated convergence. This shows that $u,v\in C^1(\mathbb{R}^n\setminus\{0\})$ and
\[
\nabla u(x)=\int_{\mathbb{R}^n}\nabla K(x-y)g(y)\,dy,\qquad \nabla K(z)=-(n-\mu)\frac{z}{|z|^2}K(z),
\]
with an absolutely convergent integral at every $x\neq0$. The same holds for $v$.

\medskip
\noindent\textit{Step 4: virial identity and localization.} Define the interaction measure $d\pi(x,y)=f(x)g(y)K(x-y)\,dx\,dy$, whose marginals are $fu\,dx$ and $gv\,dy$, and let
\[
E_f=\int\eta fu,\quad E_g=\int\eta gv,\quad E=E_f+E_g,\quad W_f=\int wfu,\quad W_g=\int wgv.
\]
Then
\begin{equation}\label{eq_virial}
a(n+\beta)E_f+b(n+\alpha)E_g-aW_f-bW_g=(n-\mu)\int_{\mathbb{R}^{2n}}\frac{(X(x)-X(y))\cdot(x-y)}{|x-y|^2}\,d\pi(x,y).
\end{equation}
We justify the identity directly, including all kernel-gradient integrals. Fix $x\ne0$, set $r=|x|$, and first suppose $r$ is sufficiently small. Since $v$ is bounded near the origin, $g(y)\le C_{u,v}r^\alpha$ on $B_{r/2}(x)$, even when $\alpha<0$. Using $|X(x)|\le Cr$, the local part satisfies
\[
|X(x)|\int_{|x-y|\le r/2}|\nabla K(x-y)|g(y)\,dy
\le C_{u,v}r^{1+\alpha}\int_0^{r/2}s^{\mu-2}\,ds
\le C_{u,v}r^{\mu+\alpha}.
\]
This is bounded as $r\downarrow0$ because $\mu>1$ and $\alpha>-\mu$. On $|x-y|>r/2$ we have $|X(x)|\le Cr\le2C|x-y|$, and therefore
\[
|X(x)|\int_{|x-y|>r/2}|\nabla K(x-y)|g(y)\,dy
\le C\int K(x-y)g(y)\,dy=Cu(x).
\]
The right side remains bounded near zero. Away from zero, the nearby source $g$ is locally bounded and $\nabla K\in L^1_{\mathrm{loc}}$, while the far part is dominated by a multiple of $K*g$. Since $X$ has compact support and $f\in L^1_{\mathrm{loc}}$, these bounds prove
\[
\int_{\mathbb R^n} f(x)|X(x)|\int_{\mathbb R^n}|\nabla K(x-y)|g(y)\,dy\,dx<\infty.
\]
Interchanging $(f,u,\beta)$ and $(g,v,\alpha)$ gives the analogous estimate. We may therefore use Fubini's theorem in both differentiated potential identities.

For $0<\rho<1$, let $\Omega_\rho=\mathbb R^n\setminus\overline B_\rho$. Integration by parts, using $f=|x|^\beta u^q$ and $\operatorname{div}(|x|^\beta X)=|x|^\beta((n+\beta)\eta-w)$, yields
\[
-\int_{\Omega_\rho}fX\cdot\nabla u
=a\int_{\Omega_\rho}((n+\beta)\eta-w)fu
+a\rho^{\beta+1}\int_{\partial B_\rho}u^{q+1}\,dS.
\]
The last term is $O_{u,v}(\rho^{n+\beta})$ and tends to zero since $\beta>-\mu>-n$. The corresponding calculation with $g$ and $v$ has a vanishing remainder $O_{u,v}(\rho^{n+\alpha})$. Taking $\rho\downarrow0$ thus gives
\[
-\int fX\cdot\nabla u-\int gX\cdot\nabla v
=a(n+\beta)E_f+b(n+\alpha)E_g-aW_f-bW_g.
\]
On the other hand, $\nabla K(z)=-(n-\mu)zK(z)/|z|^2$ and the absolute convergence proved above permit us to exchange $x$ and $y$ in the second kernel-gradient integral. The resulting sum is $(n-\mu)\int\Phi\,d\pi$, establishing \eqref{eq_virial}.

We next prove the localization estimate, following the geometric decomposition of \cite[Lemma 3.2]{OAI26}. Write
\[
\Phi(x,y)=\frac{(X(x)-X(y))\cdot(x-y)}{|x-y|^2},\qquad
D(x,y)=\frac{((x-y)\cdot x)^2}{|x-y|^2|x|^2},
\]
with $w(x)D(x,y)=0$ where $w(x)=0$. Fix $0<\delta<1/4$ and, for $d(z)>0$, set $\ell(z)=\delta d(z)^{m_0}$. Call $(x,y)$ close if $|x-y|\le\ell(z)$ for at least one endpoint $z\in\{x,y\}$ with $d(z)>0$, and far otherwise. On close pairs, since $d$ is $1$-Lipschitz, $d$ is comparable to $d(z)$ along the segment $[x,y]$. The derivative bounds for $\eta$ and $X$ in \eqref{eq_cutoff} (proved independently in Appendix \ref{app_pressure}) and the identity
\[
\nabla X=\eta I-w\,\omega_x\otimes\omega_x
\]
show by the mean value theorem and Taylor's formula that
\[
|\Phi-\eta(x)+w(x)D(x,y)|+|\eta(x)-\eta(y)|
\le C\delta d(z)^N\le C\delta(\eta(x)+\eta(y)).
\]
The marginal identities therefore bound the integral over close pairs by $C\delta E$.

For a far pair, every endpoint $z\in B_2$ satisfies $|x-y|>\ell(z)$. Using $|X|\le Cd^N$, $w\le Cd^{N-1}$, and $0\le D\le1$, we obtain
\[
|\Phi-\eta(x)+w(x)D|+|\eta(x)-\eta(y)|
\le C\sum_{z\in\{x,y\}:\,d(z)>0}
\left(d(z)^{N-1}+\frac{d(z)^N}{\ell(z)}\right).
\]
For terms associated with $x$, the tail bound \eqref{eq_tail} gives
\begin{align*}
&C\int_{B_2}f(x)\left(d(x)^{N-1}+\frac{d(x)^N}{\ell(x)}\right)
\ell(x)^{\mu-n}\,dx\\
&\quad\le C\left(\delta^{\mu-n}+\delta^{\mu-n-1}\right)
\int_{B_2}f(x)\left(d(x)^{A_1}+d(x)^{A_2}\right)\,dx\le C_\delta,
\end{align*}
where $A_1=N-1-m_0(n-\mu)=\mu n+3\mu+2n+1>0$ and
$A_2=N-m_0(n-\mu+1)=\mu n+3\mu+n-1>0$.
The terms associated with $y$ are bounded in the same way, using $g$. Finally,
$E_f-E_g=\int(\eta(x)-\eta(y))\,d\pi$ by the marginal identities.
Adding the close- and far-pair contributions and then choosing $\delta$ sufficiently small yields, for every $\varepsilon>0$,
\begin{equation}\label{eq_local}
|E_f-E_g|+\int_{\mathbb{R}^{2n}}|\Phi(x,y)-\eta(x)+w(x)D(x,y)|\,d\pi\le\varepsilon E+C_\varepsilon.
\end{equation}

\medskip
\noindent\textit{Step 5: the interval comparison.} Fix a unit vector $e$ and distinct $z,z'\in e^\perp$. This is possible because $n\ge2$. Let $\lambda=|z-z'|$ and $h=s-s'$ for $x=z+se$, $y=z'+s'e$, and put
\[
\mathcal{K}(s,s')=c_0(\lambda^2+h^2)^{(\mu-n)/2},\qquad \mathcal{D}(s,s')=\frac{h^2}{\lambda^2+h^2}.
\]
For a bounded interval $I=(s_-,s_+)$ with midpoint $m_I$, let $\mathcal{K}_I(s')=\frac12(\mathcal{K}(s_-,s')+\mathcal{K}(s_+,s'))$, and define $\mathcal{K}_J(s)$ in the second variable. The comparison \cite[Lemma 4.1]{OAI26} holds with $n-2$ replaced by $n-\mu$: for nonempty open intervals $I,J$ and coefficients $c_I,c_J\ge0$, with $I$ bounded when $c_I>0$ and $J$ bounded when $c_J>0$,
\begin{equation}\label{eq_interval}
\int_{I\times J}(c_I+c_J)\mathcal{K}\le\int_{I\times J}\Big[c_I\mathcal{K}_I+c_J\mathcal{K}_J+\big((n-\mu)\min\{c_I,c_J\}\mathcal{D}+|c_I-c_J|\big)\mathcal{K}\Big].
\end{equation}
We include the short proof, which uses only two properties of the kernel. The first is $\partial_s\mathcal{K}=-(n-\mu)h(\lambda^2+h^2)^{-1}\mathcal{K}=-\partial_{s'}\mathcal{K}$, which with integration by parts against $s-m_I$ and $s'-m_J$ gives
\[
\int_{I\times J}(2\mathcal{K}-\mathcal{K}_I-\mathcal{K}_J)=(n-\mu)\int_{I\times J}\frac{h^2-(m_I-m_J)h}{\lambda^2+h^2}\mathcal{K}.
\]
The second is that $\mathcal{K}$ depends on $h$ only through $|h|$. Let $d_0=m_I-m_J$. The image of $ds\,ds'$ under $(s,s')\mapsto h$ has density $\rho(h)=|I\cap(J+h)|$, the overlap of two intervals whose centers are $|h-d_0|$ apart. So $\rho(h)=F(|h-d_0|)$ with $F$ nonincreasing. If $d_0\ge0$ and $h>0$, then $|h-d_0|\le|h+d_0|$ and $\rho(h)\ge\rho(-h)$, with the reverse inequality when $d_0\le0$. Pairing $h$ with $-h$ against the odd function $h\mathcal{K}/(\lambda^2+h^2)$ gives $d_0\int_{I\times J}h(\lambda^2+h^2)^{-1}\mathcal{K}\ge0$, so the left side of the last display is at most $(n-\mu)\int_{I\times J}\mathcal{D}\mathcal{K}$. This proves \eqref{eq_interval} when $c_I=c_J$.

For general coefficients, if at most one of them is positive, then $|c_I-c_J|=c_I+c_J$ and \eqref{eq_interval} holds pointwise. Otherwise both intervals are bounded, all integrals are finite, and with $c=\min\{c_I,c_J\}$
\[
\int_{I\times J}\big[(c_I+c_J)\mathcal{K}-c_I\mathcal{K}_I-c_J\mathcal{K}_J\big]=c\int_{I\times J}(2\mathcal{K}-\mathcal{K}_I-\mathcal{K}_J)+(c_I-c)\int_{I\times J}(\mathcal{K}-\mathcal{K}_I)+(c_J-c)\int_{I\times J}(\mathcal{K}-\mathcal{K}_J).
\]
The first term is at most $c(n-\mu)\int\mathcal{D}\mathcal{K}$, and the last two are at most $|c_I-c_J|\int\mathcal{K}$ because $\mathcal{K}_I,\mathcal{K}_J\ge0$. When $n-\mu\le1$, the integral of $\mathcal{K}$ over an unbounded interval may diverge. This is why \eqref{eq_interval} is written with nonnegative terms only and no subtraction is made on an interval with zero coefficient.

\medskip
\noindent\textit{Step 6: the pressure estimate.} For every $\varepsilon>0$ there is $C_\varepsilon$ independent of the solution such that
\begin{equation}\label{eq_pressure}
aW_f+bW_g\le(n-\mu)\int_{\mathbb{R}^{2n}}w(x)D(x,y)\,d\pi(x,y)+\varepsilon E+C_\varepsilon.
\end{equation}
This is \cite[Proposition 5.1]{OAI26} with $n-2$ replaced by $n-\mu$. Since it is the longest step and the one most specific to the method, we give a complete proof for the kernel $K$ in Appendix \ref{app_pressure}, given \eqref{eq_tail} and \eqref{eq_interval}.

\medskip
\noindent\textit{Step 7: conclusion.} Write $T=\int w(x)D(x,y)\,d\pi$. Inserting \eqref{eq_local} and $\int\eta(x)\,d\pi=E_f$ into \eqref{eq_virial} gives
\[
c_1E_f+c_2E_g\le aW_f+bW_g-(n-\mu)T+(n-\mu)(\varepsilon E+C_\varepsilon),\qquad c_1=a(n+\beta)-(n-\mu),\quad c_2=b(n+\alpha).
\]
By \eqref{eq_pressure}, the right side is at most $(n-\mu+1)\varepsilon E+C_\varepsilon$. Since $c_1+c_2=\gamma_\mu$, we have $c_1E_f+c_2E_g=\frac{\gamma_\mu}2E+\frac{c_1-c_2}2(E_f-E_g)$, and \eqref{eq_local} yields
\[
\frac{\gamma_\mu}{2}E\le\Big(n-\mu+1+\frac{|c_1-c_2|}{2}\Big)\varepsilon E+C_\varepsilon.
\]
The energy $E$ is finite for each solution, because $u,v$ are bounded on $B_2$ and $\alpha,\beta>-\mu>-n$. Choosing $\varepsilon$ small, we obtain $E\le C$ with $C$ independent of the solution. As $\eta=1$ on $B_1$,
\[
\int_{B_1}(fu+gv)\le C
\]
for every positive continuous solution of \eqref{eq_int}. Applying this to $(u_R,v_R)$ and using \eqref{eq_energy_scaling} with $\sigma_1+\sigma_2+\mu-n>0$, we get $\int_{B_R}(fu+gv)\le CR^{n-\mu-\sigma_1-\sigma_2}\to0$ as $R\to\infty$. But $\int_{B_R}(fu+gv)\ge\int_{B_1}(fu+gv)>0$ for $R\ge1$. This contradiction proves Theorem \ref{thm_integral}.

\section{Positive solutions without the super polyharmonic property}\label{sec_example}

\begin{proof}[Proof of Proposition \ref{prop_example}]
Since $m$ is even, $(-\Delta)^m=\Delta^m$. We look for radial functions $w_i=\Delta^iu$ and $z_i=\Delta^iv$, $0\le i\le m-1$, solving
\begin{equation}\label{eq_ode}
w_i(r)=1+\int_0^rt^{1-n}\int_0^t\tau^{n-1}w_{i+1}(\tau)\,d\tau\,dt,\qquad z_i(r)=1+\int_0^rt^{1-n}\int_0^t\tau^{n-1}z_{i+1}(\tau)\,d\tau\,dt,
\end{equation}
with $w_m=z_0^p$ and $z_m=w_0^q$. For a continuous function $h$ on $[0,\infty)$, the function $W(r)=1+\int_0^rt^{1-n}\int_0^t\tau^{n-1}h(\tau)\,d\tau\,dt$ satisfies $W(0)=1$ and $W'(0)=0$, and $W(|x|)$ solves $\Delta W=h$ in the sense of distributions in $\mathbb{R}^n$. Hence a continuous solution of \eqref{eq_ode} on $[0,\rho)$ with $w_0,z_0\ge1$ gives $u=w_0$, $v=z_0$ with $\Delta^mu=v^p$ and $\Delta^mv=u^q$ in $B_\rho$. Since $u,v\ge1$, the nonlinearities are smooth functions of $u$ and $v$, and Schauder estimates give $u,v\in C^\infty(B_\rho)$.

The right sides of \eqref{eq_ode} are Lipschitz in the unknowns on the set where $w_0,z_0\ge\frac12$. The contraction principle gives a solution on some $[0,\rho)$, and it can be continued as long as all $w_i,z_i$ stay bounded and $w_0,z_0$ stay away from zero. On the existence interval, $w_m,z_m>0$, so by \eqref{eq_ode} the functions $w_{m-1}$, $z_{m-1}$ are nondecreasing and at least $1$. Descending, every $w_i$ and $z_i$ is nondecreasing and at least $1$. In particular $u,v\ge1$ and $\Delta u=w_1\ge1$, $\Delta v=z_1\ge1$.

To exclude a finite maximal radius $\rho$, we bound the solution on $[0,\rho)$. Since $z_0^p$ is nondecreasing, $\int_0^t\tau^{n-1}z_0^p\le t^nz_0(r)^p/n$ for $t\le r$, and iterating in \eqref{eq_ode} gives
\[
w_0(r)\le(1+r^2)^m+r^{2m}z_0(r)^p,\qquad z_0(r)\le(1+r^2)^m+r^{2m}w_0(r)^q\qquad(0\le r<\rho).
\]
If $pq<1$, then with $P=(1+\rho^2)^m$ we obtain $w_0\le P+2^p\rho^{2m}(P^p+\rho^{2mp}w_0^{pq})$ on $[0,\rho)$. This bounds $w_0$, hence $z_0$, hence every $w_i$ and $z_i$ on $[0,\rho)$, which is impossible. If $p=q=1$, then \eqref{eq_ode} is a linear Volterra system with continuous kernels, and its solution exists on $[0,\infty)$. In both cases the solution is global.

Finally, $pq\le1$ gives $1-\frac1{p+1}-\frac1{q+1}=\frac{pq-1}{(p+1)(q+1)}\le0$. Hence $\frac{n}{p+1}+\frac{n}{q+1}\ge n>n-2m$, which is \eqref{eq_sub}.
\end{proof}

\appendix

\section{Proof of the pressure estimate}\label{app_pressure}

In this appendix we prove \eqref{eq_pressure}. The argument is that of \cite[Section 5]{OAI26}, written for the kernel $K(z)=c_0|z|^{\mu-n}$. We keep the notation of Section \ref{sec_integral}. In particular $(u,v)$ is a positive continuous solution of \eqref{eq_int} with $n\ge2$ and $1<\mu<n$, and $f=|x|^\beta u^q$, $g=|x|^\alpha v^p$, with the understanding that, in layer-set constructions only, the source functions have been redefined to vanish at the single point $0$. Every constant $C$ below depends only on $n,\mu,c_0,p,q,\alpha,\beta$ and on the parameters indicated as subscripts. It never depends on the solution or on the direction $e$.

\medskip
\noindent\textit{The cutoff.} Recall the cutoffs $d$, $\eta$, $w$, $X$ and the integers $m_0=n+3$, $k=n(m_0+1)$, $N=k+n+2$ from Section \ref{sec_integral}.
On $1<|x|<2$ we have $d=1-s$ with $s=|x|-1$, and $1-s^3=d(1+s+s^2)$. Differentiating gives
\begin{equation}\label{eq_cutoff}
d^N\le\eta\le Cd^N,\qquad|\nabla\eta|+w\le Cd^{N-1},\qquad|\nabla^2\eta|+|\nabla w|+|\nabla^2X|\le Cd^{N-2},\qquad dw\le C\eta,
\end{equation}
while $\eta=1$ and $w=0$ on $B_1$, and $\eta=w=0$ outside $B_2$. Since $w$ vanishes to second order at $|x|=1$, the function $w$ is $C^1$ on $\mathbb{R}^n$.

\medskip
\noindent\textit{Cap weights.} Fix $\kappa\in(0,1)$ and $\delta\in(0,1/4)$, and put $\ell(x)=\delta d(x)^{m_0}$ when $d(x)>0$. Let $\mathbb{E}_e$ and $\mathbb{P}_e$ denote expectation and probability for the normalized surface measure on the unit sphere $\mathbb{S}^{n-1}$. Fix $\psi\in C^\infty([0,\infty))$ with $0\le\psi\le1$, $\psi=1$ near $0$, and $\psi=0$ on $[1,\infty)$, and set $Z(l)=\mathbb{E}_e\psi(|e-\theta|/l)$ for $0<l\le1$. By rotation invariance, $Z$ does not depend on $\theta\in\mathbb{S}^{n-1}$, and polar coordinates around $\theta$ give $cl^{n-1}\le Z(l)\le Cl^{n-1}$ and $|Z'(l)|\le Cl^{n-2}$. For $e\in\mathbb{S}^{n-1}$ define
\[
h_e(x)=w(x)\frac{\psi\big(|e-\omega_x|/(\kappa d(x))\big)}{Z(\kappa d(x))}\quad(1<|x|<2),\qquad h_e(x)=0\quad\text{otherwise},
\]
where $\omega_x=x/|x|$. Then $\mathbb{E}_eh_e=w$. Differentiating the cap profile or $Z(\kappa d)^{-1}$ costs at most a relative factor $C_\kappa d^{-1}$, so by \eqref{eq_cutoff}
\begin{equation}\label{eq_cap}
0\le h_e\le C_\kappa d^{N-n},\qquad|\nabla h_e|\le C_\kappa d^{N-n-1}\qquad\text{on }1<|x|<2.
\end{equation}
The profile is constant near its center, so $h_e$ is $C^1$ there. Moreover, $h_e$ tends to zero at both boundary spheres, because $w$ vanishes at $|x|=1$ and $N-n>1$. Hence $h_e$ is continuous on $\mathbb{R}^n$, and $|\nabla h_e|\le C_\kappa d^{N-n-1}$ at every point off the spheres $|x|=1$ and $|x|=2$, since $h_e=0$ outside the annulus. A segment meets these spheres in at most four points, so this bound may be integrated along segments.

\begin{lemma}\label{lem_angular}
For every direction $e$,
\begin{equation}\label{eq_ang1}
\int_{\mathbb{R}^{2n}}|h_e(x)-h_e(y)|\,d\pi(x,y)\le C_\kappa\delta E+C_{\kappa,\delta}.
\end{equation}
Moreover, with $D_e(x,y)=((x-y)\cdot e)^2/|x-y|^2$,
\begin{equation}\label{eq_ang2}
\mathbb{E}_e\int_{\mathbb{R}^{2n}}h_e(x)D_e(x,y)\,d\pi(x,y)\le\int_{\mathbb{R}^{2n}}w(x)D(x,y)\,d\pi(x,y)+C\kappa E,
\end{equation}
with $C$ independent of $\kappa$ and $\delta$.
\end{lemma}

\begin{proof}
Call a pair $(x,y)$ close if $|x-y|\le\ell(z)$ for some endpoint $z\in\{x,y\}$ with $d(z)>0$, and far otherwise. Since $d$ is $1$-Lipschitz and $\ell(z)\le\delta d(z)$, for a close pair we have $(1-\delta)d(z)\le d\le(1+\delta)d(z)$ on the segment $[x,y]$. Integrating the gradient bound along the segment gives
\[
|h_e(x)-h_e(y)|\le C_\kappa d(z)^{N-n-1}\ell(z)=C_\kappa\delta d(z)^{N+2}\le C_\kappa\delta\big(\eta(x)+\eta(y)\big),
\]
since $m_0=n+3$ and $d(z)^N\le\eta(z)$. By the marginal identities, close pairs contribute at most $C_\kappa\delta(E_f+E_g)$. For a far pair we use $|h_e(x)-h_e(y)|\le h_e(x)+h_e(y)$. The term $h_e(x)$ vanishes unless $1<|x|<2$, and then $|x-y|>\ell(x)$. By \eqref{eq_cap} and \eqref{eq_tail}, its contribution is at most
\[
C_\kappa\int_{B_2}f(x)d(x)^{N-n}\int_{|x-y|>\ell(x)}K(x-y)g(y)\,dy\,dx\le C_\kappa\delta^{\mu-n}\int_{B_2}f(x)d(x)^{N-n-m_0(n-\mu)}\,dx\le C_{\kappa,\delta},
\]
because $N-n-m_0(n-\mu)=\mu n+3\mu+n+2>0$ and $\int_{B_2}f\le C$. The term $h_e(y)$ is treated in the same way with $g$ in place of $f$. This proves \eqref{eq_ang1}.

For unit vectors $\hat z$, $e$, $\omega$ we have $|(\hat z\cdot e)^2-(\hat z\cdot\omega)^2|\le|e-\omega||e+\omega|\le2|e-\omega|$. On the support of $h_e(x)$ we have $|e-\omega_x|\le\kappa d(x)$, so $h_eD_e\le h_eD+2\kappa dh_e$. Averaging in $e$, using $\mathbb{E}_eh_e=w$ and $dw\le C\eta$, gives \eqref{eq_ang2}.
\end{proof}

\medskip
\noindent\textit{Layer intervals.} Fix $e$. For $x\in\mathbb{R}^n$ and $0<t<f(x)$, let $I=I(e,x,t)$ be the connected component containing $x$ of the open set $\{x+se:\ f(x+se)>t\}$ on the line $x+\mathbb{R}e$. For the purpose of defining the layers, we assign the value $f(0)=0$ (and similarly $g(0)=0$), even if a weight exponent is zero. This modification at one point changes none of the integrals or potentials. Since $f$ is continuous away from the origin, the resulting layer set is open. Define
\begin{equation}\label{eq_chi}
\chi_f(e,x,t)=\begin{cases}\inf_Ih_e,&\text{if $I$ is bounded and $|I|\le2\delta\inf_Id$},\\0,&\text{otherwise},\end{cases}
\end{equation}
and define $\chi_g(e,y,t')$ in the same way from the layers $\{g>t'\}$. Thus $\chi_f(e,x,t)$ depends on $x$ only through $I$, and $0\le\chi_f(e,x,t)\le h_e(x)$.

These functions are jointly measurable in $(e,x,t)$. Indeed, the forward exit distance of $x$ from its component exceeds $r\ge0$ exactly when the compact segment $\{x+se:0\le s\le r\}$ lies in $\{f>t\}$, and this is an open condition in $(e,x,t)$. The same holds for the backward exit distance. The length of $I$ is their sum, possibly $+\infty$. Since $d$ and $h_e$ are continuous, the infimum in \eqref{eq_chi} may be taken over rational points between the forward and backward exit points. When $I$ is bounded, its endpoints are $x$ plus or minus the exit distances times $e$, so they are measurable as well, and so is the endpoint average $u_I$ defined below, because $u$ is continuous. Define the losses
\[
L_f(e)=\int_{\mathbb{R}^n}u(x)\int_0^{f(x)}\big(h_e(x)-\chi_f(e,x,t)\big)\,dt\,dx,\qquad L_g(e)=\int_{\mathbb{R}^n}v(y)\int_0^{g(y)}\big(h_e(y)-\chi_g(e,y,t')\big)\,dt'\,dy.
\]
They are finite because $h_e$ is bounded and supported in $\overline{B_2}\setminus B_1$, where $u$ and $v$ are bounded and $f,g$ are integrable.

\begin{lemma}\label{lem_loss}
For every $P\ge1$,
\[
\mathbb{E}_e\big(L_f(e)+L_g(e)\big)\le\big(C_\kappa\delta+C_{\kappa,\delta}P^{-1}\big)E+CP,
\]
where the last constant does not depend on $\kappa$, $\delta$, or $P$.
\end{lemma}

\begin{proof}
We treat $L_f$. Only points with $1<|x|<2$ contribute. Fix such an $x$ and write $d=d(x)$ and $\ell=\delta d^{m_0}$.

If $|I|\le\ell$, then $\inf_Id\ge d-\ell\ge(1-\delta)d$ and $|I|\le\delta d\le2\delta\inf_Id$, so the component is retained. Integrating \eqref{eq_cap} along $I$ as in Lemma \ref{lem_angular} gives $h_e(x)-\chi_f(e,x,t)\le C_\kappa\delta d^N$.

If $|I|>\ell$, then one of the two segments of length $\ell/2$ starting at $x$ in the directions $\pm e$ lies in $\{f>t\}$. Let $F$ be the set of $e$ for which the forward segment does. Polar coordinates centered at $x$ give
\[
\int_{B_{\ell/2}(x)}f\ge\frac{|\mathbb{S}^{n-1}|}{n}\,t\Big(\frac{\ell}{2}\Big)^n\,\mathbb{P}_e(F).
\]
Since $B_{\ell/2}(x)\subset B_4$, \eqref{eq_tail} gives $\mathbb{P}_e\{|I|>\ell\}\le2\mathbb{P}_e(F)\le C/(t\ell^n)$. For $t\ge Pd^{-k}$ we bound the loss on this event by $\sup_eh_e(x)\le C_\kappa d^{N-n}$, and obtain
\[
\mathbb{E}_e\big[h_e(x)-\chi_f(e,x,t)\big]\le C_\kappa\delta d^N+C_\kappa d^{N-n}\frac{C}{t\delta^nd^{m_0n}}\le\big(C_\kappa\delta+C_{\kappa,\delta}P^{-1}\big)d^N,
\]
because $k-m_0n=n$. For $t<Pd^{-k}$ we use only $\mathbb{E}_e\big[h_e(x)-\chi_f\big]\le\mathbb{E}_eh_e(x)=w(x)$, so
\[
\int_0^{\min\{f(x),Pd^{-k}\}}\mathbb{E}_e\big[h_e(x)-\chi_f(e,x,t)\big]\,dt\le Pd^{-k}w(x)\le CPd^{N-1-k}=CPd^{n+1}\le CP.
\]
By Tonelli's theorem we may integrate these bounds in $t$ and then against $u(x)\,dx$. The first gives at most $(C_\kappa\delta+C_{\kappa,\delta}P^{-1})\int d^Nfu\le(C_\kappa\delta+C_{\kappa,\delta}P^{-1})E_f$, and the second at most $CP\int_{B_4}u\le CP$. The same argument applies to $L_g$.
\end{proof}

\medskip
\noindent\textit{Endpoint values.} Suppose $\chi_f(e,x,t)>0$, and let $x_\pm$ be the endpoints of $I$. Then $I$ is bounded and $h_e\ge\chi_f>0$ on $I$. Hence, by continuity, $h_e(x_\pm)\ge\chi_f>0$, so $x_\pm$ lie in the open annulus $1<|x|<2$. There $f$ is continuous, and maximality of $I$ gives $f(x_\pm)=t$, that is, $u(x_\pm)=(t|x_\pm|^{-\beta})^{1/q}$. Put $u_I=\frac12\big(u(x_-)+u(x_+)\big)$ and $\tau=(t|x|^{-\beta})^{1/q}$. Since $|x_\pm-x|\le|I|\le2\delta d(x)$ and the function $r\mapsto r^{-\beta/q}$ has a bounded logarithmic derivative on $[1,2]$,
\begin{equation}\label{eq_endpoint}
|u_I-\tau|\le C\delta d(x)\tau\le C\delta d(x)u(x),
\end{equation}
where we used $t<f(x)$, that is, $\tau<u(x)$. The same holds for $g$, with $v_J$ and $(t'|y|^{-\alpha})^{1/p}$. When a component weight vanishes, every product involving its endpoint average is defined to be zero.

\medskip
\noindent\textit{Assembly over layers.} Fix $e$ and $t,t'>0$. Write $x=z+se$ and $y=z'+s'e$ with $z,z'\in e^\perp$, and take $z\neq z'$. For each component $I$ of the $f$-layer on the line through $z$ and each component $J$ of the $g$-layer on the line through $z'$, apply \eqref{eq_interval} with $c_I=\chi_I$ and $c_J=\chi_J$, the common values of $\chi_f$ on $I$ and of $\chi_g$ on $J$. Since $\min\{\chi_I,\chi_J\}\le\chi_I\le h_e(x)$ for $x\in I$,
\begin{equation}\label{eq_pair}
\int_{I\times J}(\chi_I+\chi_J)\mathcal{K}\le\int_{I\times J}\Big[\chi_I\mathcal{K}_I+\chi_J\mathcal{K}_J+\big((n-\mu)h_e(x)D_e(x,y)+|\chi_I-\chi_J|\big)\mathcal{K}\Big]\,ds\,ds'.
\end{equation}
Every term is nonnegative. Each open layer on a line has at most countably many components. We sum \eqref{eq_pair} over all pairs of components and integrate over $z,z'\in e^\perp$ and $t,t'>0$. The set $\{z=z'\}$ is null because $n\ge2$. By Tonelli's theorem and $\int_0^\infty\mathbf{1}_{\{g(y)>t'\}}dt'=g(y)$, the identity $u=K*g$ turns the two terms on the left into
\[
U_f(e)=\int_{\mathbb{R}^n}\int_0^{f(x)}\chi_f(e,x,t)u(x)\,dt\,dx,\qquad U_g(e)=\int_{\mathbb{R}^n}\int_0^{g(y)}\chi_g(e,y,t')v(y)\,dt'\,dy.
\]
For the endpoint terms, integrating $\mathcal{K}_I$ over all $g$-layers recovers $\frac12\big((K*g)(x_-)+(K*g)(x_+)\big)=u_I$. For each endpoint, the omitted points $y$ form the line $z+\mathbb{R}e$, which is a null set. Hence these terms become
\[
V_f(e)=\int_{\mathbb{R}^n}\int_0^{f(x)}\chi_f(e,x,t)u_I\,dt\,dx,\qquad V_g(e)=\int_{\mathbb{R}^n}\int_0^{g(y)}\chi_g(e,y,t')v_J\,dt'\,dy.
\]
The pressure term becomes $(n-\mu)\int h_e(x)D_e(x,y)\,d\pi$. We check finiteness before subtracting anything. Let $H_f(e)=\int h_efu$ and $H_g(e)=\int h_egv$. Both are finite for each solution. We have $U_f\le H_f$, $U_g\le H_g$, and by \eqref{eq_endpoint} also $V_f\le CH_f$ and $V_g\le CH_g$. The pressure term is at most $(n-\mu)H_f$. Finally, $|\chi_f(e,x,t)-\chi_g(e,y,t')|\le h_e(x)+h_e(y)$, so the integrated mismatch term is at most $H_f+H_g$. All integrated terms are therefore finite, even when $\int_J\mathcal{K}=\infty$ for some unbounded component. Subtracting $V_f+V_g$ and using
\[
|\chi_f(e,x,t)-\chi_g(e,y,t')|\le|h_e(x)-h_e(y)|+\big(h_e(x)-\chi_f(e,x,t)\big)+\big(h_e(y)-\chi_g(e,y,t')\big),
\]
we find that $S_e=U_f(e)-V_f(e)+U_g(e)-V_g(e)$ satisfies
\begin{equation}\label{eq_Se}
S_e\le(n-\mu)\int h_e(x)D_e(x,y)\,d\pi+\int|h_e(x)-h_e(y)|\,d\pi+L_f(e)+L_g(e).
\end{equation}

\medskip
\noindent\textit{Recovery of the pressure.} For $x\neq0$,
\[
\int_0^{f(x)}\Big[u(x)-\big(t|x|^{-\beta}\big)^{1/q}\Big]\,dt=f(x)u(x)-\frac{q}{q+1}|x|^{-\beta/q}f(x)^{1+1/q}=af(x)u(x),
\]
since $|x|^{-\beta/q}f(x)^{1/q}=u(x)$. With $\tau=(t|x|^{-\beta})^{1/q}\le u(x)$ and \eqref{eq_endpoint},
\[
h_e(x)(u-\tau)-\chi_f(u-u_I)=(h_e(x)-\chi_f)(u-\tau)+\chi_f(u_I-\tau)\le(h_e(x)-\chi_f)u+C\delta d(x)h_e(x)u,
\]
where $u=u(x)$ and $\chi_f=\chi_f(e,x,t)$. Integrating over $0<t<f(x)$ and then in $x$, and doing the same for $g$, we obtain
\[
aH_f(e)+bH_g(e)\le S_e+L_f(e)+L_g(e)+C\delta\int d\,h_e\,(fu+gv)\,dx.
\]
We insert \eqref{eq_Se}, which contains a second copy of $L_f+L_g$, and average over $e$. Since $\mathbb{E}_eH_f=W_f$, $\mathbb{E}_eH_g=W_g$, and $\mathbb{E}_e(dh_e)=dw\le C\eta$, Lemmas \ref{lem_angular} and \ref{lem_loss} give
\[
aW_f+bW_g\le(n-\mu)\int w(x)D(x,y)\,d\pi+\big(C\kappa+C_\kappa\delta+C_{\kappa,\delta}P^{-1}+C\delta\big)E+C_{\kappa,\delta}+CP.
\]
The constant multiplying $\kappa$ does not depend on $\kappa$, $\delta$, or $P$. The one multiplying $\delta$ depends only on $\kappa$, and the one multiplying $P^{-1}$ does not depend on $P$. We first choose $\kappa$, then $\delta$, then $P$, so that the coefficient of $E$ is at most $\varepsilon$. This proves \eqref{eq_pressure}.

\section*{Declarations}

\smallskip
\noindent\textbf{Conflict of interest.} The author declares that he has no conflict of interest.

\smallskip
\noindent\textbf{Data availability.} Data sharing is not applicable to this article, as no datasets were generated or analyzed during the current study.

\smallskip
\noindent\textbf{Use of AI tools.} AI-assisted tools were used during the preparation of the manuscript. The author reviewed and edited all AI-assisted output and takes full responsibility for the content of the work.

\end{document}